\documentclass[11pt]{article} 
\title{{\bf Smooth Rigidity of Uniformly Quasiconformal Anosov Flows}} 
\author{Yong Fang} 
\date{{\it Laboratoire de Math\'ematique d'Orsay, U.M.R. 
8628 du C.N.R.S, Universit\'e Paris-sud, France} \\ {\it (e-mail: 
fangyong1@yahoo.fr)}} 

\usepackage{amsmath,amsthm,amsfonts,amssymb,amscd,amstext} 
\chardef\bslash=`\\

\theoremstyle{definition}

\theoremstyle{remark}

\begin{document} 
\maketitle 
\renewcommand{\sectionmark}[1]{} 
{\bf Abstract}--{\it We classify the $C^\infty$ volume-preserving uniformly quasiconformal 
Anosov flows, such that $E^+\oplus E^-$ is $C^\infty$ and the dimensions of $E^+$ and $E^-$ are 
at least two. Then we deduce a classification of volume-preserving uniformly quasiconformal 
Anosov flows with smooth distributions.}\\\\
{\bf 1.}  {\bf Introduction}

Conformal geometry is a classically and currently fascinating subject, which is meaningfully mixed 
together with hyperbolic dynamical systems under the impulsion of the classical {\bf [Su]}, {\bf [Ka]} 
and {\bf [Yu]} and the recent {\bf [Sa]} and {\bf [K-Sa]}. In this paper, we study the rigidity of such 
systems. 

Let $M$ be a $C^\infty$-closed manifold. A $C^\infty$-flow, $\phi_t$, generated by a non-singular vector field $X$ 
on $M$ is called 
Anosov, if there exists a $\phi_t$-invariant splitting of the tangent bundle
$$TM =\mathbb{R}X\oplus E^+\oplus E^-,$$
a Riemannian metric on $M$ and two positive numbers $a$ and $b$, such that 
$$\forall \  u^\pm\in E^\pm,\  \forall \  t > 0, \  \parallel D\phi_{\mp t}(u^\pm)\parallel\leq 
a\cdot e^{-bt}\parallel u^\pm\parallel.$$

Since $M$ is compact, then the definition of an Anosov flow is independent of the Riemannian metric 
chosen. A metric is called Lyapunov, if with respect to it, the constant $a$ above can be taken to be $1$. 
Then for any Anosov flow, there exists always a Lyapunov metric (see {\bf [Sh]}). $E^+$ and $E^-$ are 
called the unstable and stable distributions of $\phi_t$. The 
{\it canonical $1$-form} of $\phi_t$ 
is by definition the continuous $\phi_t$-invariant $1$-form $\lambda$ such that 
$$\lambda(E^+\oplus E^-) \equiv 0,\  \  \lambda(X) \equiv 1.$$
Denote $E^\pm\oplus \mathbb{R}X$ 
by $E^{\pm, 0}$. Then $E^{+, 0}$ and $E^{-, 0}$ are called the weak unstable and weak stable 
distributions of $\phi_t$. $E^\pm$ and $E^{\pm, 0}$ are 
all integrable to continuous foliations with $C^\infty$ leaves, denoted respectively by 
$\mathcal{F}^\pm$ and $\mathcal{F}^{\pm, 0}$ (see {\bf [HK]}). The corresponding leaves, passing through $x$, 
are denoted by $W^\pm_x$ and $W^{\pm, 0}_x.$ 

Define two functions on $M\times\mathbb{R}$ as follows,
$$K^\pm(x, t) := \frac{ max \{ \parallel D\phi_t(u)\parallel \  \mid u\in E^\pm_x,\  \parallel u\parallel =1\}}
{min \{ \parallel D\phi_t(u)\parallel \  \mid u\in E^\pm_x, \  \parallel u\parallel =1\}}.$$
If $K^-$ $(K^+)$ is bounded, then the Anosov flow $\phi_t$ is called uniformly 
quasiconformal on the stable (unstable) distribution. If $K^+$ and $K^-$ 
are both bounded, then $\phi_t$ is called uniformly quasiconformal. 
The corresponding notions for Anosov diffeomorphisms are defined similarly (see {\bf [Sa]}).

Among such uniformly quasiconformal Anosov systems, the uniformly quasiconformal 
geodesic flows were classically studied. In {\bf [Ka]} and {\bf [Yu]}, some elegant 
rigidity results were obtained via the sphere at infinity. Quite recently in {\bf [Sa]}, V. Sadovskaya 
obtained a classification of uniformly quasiconformal contact Anosov flows of dimension at least 
five. In this paper, we improve all these rigidity results by proving  

$\ $

{\bf Theorem 1.} {\it Let $\phi_t$ be a $C^\infty$ volume-preserving 
uniformly quasiconformal Anosov flow on a closed manifold $M$. If 
$E^+\oplus E^-$ is $C^\infty$ and the dimensions of $E^+$ and $E^-$ are at least 2, then 
up to a constant change of time scale 
and finite covers, $\phi_t$ is $C^\infty$ flow equivalent either to the suspension of 
a hyperbolic automorphism of a torus, or to a canonical perturbation of the geodesic flow 
of a hyperbolic manifold.}

$\ $

For each vector field $X$ on $M$, a canonical perturbation 
of the flow of $X$ is by definition the flow of $\frac{X}{1 +\alpha(X)}$, where 
$\alpha$ is a $C^\infty$ closed $1$-form on $M$ such that $1+\alpha(X) >0$. 

To prove Theorem $1$, 
our main technique is the following Theorem $2.$ Let us recall at first some terminology. 
A flow, $\psi_t$, is called topologically transitive, if it has a dense orbit. It is called 
topologically mixing, if for $\forall \  U, \  V$, nonempty open subsets of $M$, $\exists \  T >0,$ such that 
$$\psi_t U\cap V\not =\emptyset,\forall\  t\geq T.$$ 

Among the dynamical invariants of an 
Anosov flow $\phi_t$, its topological entropy is the most fundamental one, denoted by 
$h_{top}(\phi_t)$. Denote by $\mathcal{M}(\phi_t)$ the set of $\phi_t$-invariant probability measures. For $\forall\  
\nu\in \mathcal{M}(\phi_t)$, denote by $h_{\nu}(\phi_t)$ the metrical entropy of $(\phi_t,\  \nu)$. 
It is classically shown (see {\bf [HK]}) that if $\phi_t$ is topologically transitive, then 
there exists a unique $\phi_t$-invariant probability measure, $\mu$, such that 
$$h_{\mu}(\phi_t) =h_{top}(\phi_t) =sup_{\nu\in \mathcal{M}(\phi_t)}\{ h_{\nu}(\phi_t) \}.$$
This measure $\mu$ is called the Bowen-Margulis measure of $\phi_t$.

$\ $

{\bf Theorem 2.} {\it Let $\phi_t$ be a $C^\infty$ topologically transitive 
uniformly quasiconformal Anosov flow on 
a closed manifold, such that $E^+\oplus E^-$ is $C^\infty$ and the dimensions of $E^+$ and $E^-$ 
are at least 2. If its Bowen-Margulis measure is in the Lebesgue measure class, then, up to a constant change 
of time scale and finite covers, $\phi_t$ is $C^\infty$ flow equivalent either to the suspension 
of a hyperbolic automorphism of a torus, or to the geodesic flow of a hyperbolic manifold.}

$\ $

By considering the suspensions, we can draw from Theorem $1.$ the following

$\ $

{\bf Corollary 1.} {\it Let $\phi$ be a $C^\infty$ volume-preserving uniformly quasiconformal Anosov 
diffeomorphism on a closed manifold 
$\Sigma$. If the dimensions 
of $E^+$ and $E^-$ are at least $2$, then up to finite covers, 
$\phi$ is $C^\infty$ conjugate to a hyperbolic automorphism of a torus.}

$\ $

In {\bf [K-Sa]}, a similar proposition is proved for the topologically transitive case 
by assuming that the dimensions of $E^\pm$ are at least three.

Recently, P. Foulon proved an entropy rigidity theorem for three dimensional 
contact Anosov flows (see {\bf [Fo]}). In the case of dimension three, an 
Anosov flow is certainly uniformly quasiconformal. So by combining his result 
with our Theorem $2$, we obtain the following 

$\ $

{\bf Corollary 2.} {\it Let $\phi_t$ be a $C^\infty$ uniformly quasiconformal contact Anosov flow. 
If its Bowen-Margulis measure is in the Lebesgue class, then up to 
a constant change of time scale and finite covers, $\phi_t$ is $C^\infty$ flow equivalent to the 
geodesic flow of a hyperbolic manifold.}

$\ $

By extending partially our Theorem $1.$ to the case of codimension one, we get the following theorem 
which generalizes the classification result in {\bf [Gh]}.   

$\ $

{\bf Theorem 3.} {\it Let $\phi_t$ be a $C^\infty$ volume-preserving uniformly quasiconformal 
Anosov flow. If $E^+$ and $E^-$ are both $C^\infty$, then up to a constant change of time scale 
and finite covers, $\phi_t$ is $C^\infty$ flow equivalent either to the suspension of 
a hyperbolic automorphism of a torus, or to a canonical perturbation of the geodesic flow 
of a hyperbolic manifold.}

$\ $

In particular, we deduce the following 

$\ $

{\bf Corollary 3.} {\it Let $\phi$ be a $C^\infty$ volume-preserving uniformly quasiconformal Anosov 
diffeomorphism. If $E^+$ and $E^-$ are both $C^\infty$, then up to finite covers, $\phi$ is $C^\infty$ 
conjugate to a hyperbolic automorphism of a torus.}

$\ $

In Section $2.$ below, we fix the terminology and we recall and prove several dynamical 
and geometric lemmas. 
In Section $3$, we prove Theorem $2.$ Then in Section $4$, we prove Theorem $1.$ by reducing it 
to the case of Theorem $2$. In 
addition, we deduce Corollary $1.$ and prove Theorem $3$.\\\\
{\bf 2. Preliminaries.}\\
{\bf 2.1. Linearizations and smooth conformal structures.}

In this subsection, we review and adapt some results of {\bf [Sa]} to our situation. The starting point is 
the following elegant proposition from {\bf [Sa]}, which generalizes a one-dimensional result of {\bf [KL]}.

$\ $

{\bf Proposition 2.1.1.} {\it Let $f$ be a diffeomorphism of a compact Riemannian manifold $M$, and let $W$ 
be a continuous f-invariant foliation with $C^\infty$ leaves. Suppose that $\parallel Df\mid_{TW}\parallel <1,$ 
and there exist $C>0$ and $\epsilon >0$ such that for any $x\in M$ and $n\in \mathbb{N}$, 
$$\parallel (Df^n\mid_{T_xW})^{-1}\parallel\cdot\parallel Df^n\mid_{T_xW}\parallel^2\  \leq C(1-\epsilon)^n.$$
Then for any $x\in M$, there exists a $C^\infty$ diffeomorphism $h_x : W_x\to T_xW$ such that\\ 
$(i).$ $h_{fx}\circ f=Df_x\circ h_x,$\\
$(ii).$ $ h_x(x) =0$ and $(Dh_x)_x$ is the identity map,\\ 
$(iii).$ $h_x$ depends continuously on $x$ in $C^\infty$ topology.\\
In addition, the family $h$ of maps $h_x$ satisfying $(i)$, $(ii)$ and $(iii)$ is unique.}

$\ $          
 
Now let $\phi_t$ be a uniformly quasiconformal Anosov flow on a closed manifold $M$. Fix a Lyapunov 
metric on $M$. Then for $\forall\  s <0,$ 
$$\parallel D\phi_s\mid_{E^+}\parallel <1.$$
 
Since the flow is uniformly quasiconformal, then for each 
negative number $s$, $(\phi_s,\  \mathcal{F}^+)$ satisfies the 
conditions of the previous proposition. So for $\forall \  x\in M$, 
there exists a $C^\infty$ diffeomorphism $h^{+, s}_x : W^+_x\rightarrow E^+_x$, such that \\
$(i)$. $h^{+, s}_{\phi_s(x)}\circ \phi_s = D_x\phi_s\circ h^{+, s}_x,$\\
$(ii)$. $h^{+, s}_x(x) =0$ and $(Dh^{+, s}_x)_x$ is the identity map,\\
$(iii)$. $h^{+, s}_x$ depends continuously on $x$ in the $C^\infty$ topology.

For $\forall \  m\in \mathbb{N}$, we observe easily that 
$\{ h^{+, \frac{s}{m}} \}_{x\in M}$ satisfies also these three conditions with respect to $\phi_s$. Then by the 
uniqueness 
of this family of maps for $\phi_s$, we get 
$$h^{+, \frac{s}{m}}_x =h^{+, s}_x,\  \forall \  x\in M.$$

We deduce that for $\forall \  a\in \mathbb{Q}$ and $a <0$, $h^{+, a}_x =h^{+, -1}_x, \  \forall\  x\in M.$ Then 
by the condition $(iii)$, we get 
$$ h^{+, -1}_{\phi_t(x)}\circ \phi_t =(D_x\phi_t)\circ h^{+, -1}_x, \  \forall \  x\in M, \  \forall \  t <0.$$
Denote $h^{+, -1}_x$ by $h^+_x$. Thus we have 
$$ h^+_{\phi_t(x)}\circ \phi_t=(D_x\phi_t)\circ h^+_x,\  \forall\  x\in M,\  \forall \  t\in \mathbb{R}.$$ 
This continuous family of $C^\infty$ maps $\{ h^+_x\}_{x\in M}$ is called 
the {\it unstable linearization} of $E^+.$ Similarly, 
we get the {\it stable linearization} of $E^-$, $\{ h^-_x \}_{x\in M}$. 

For the sake of completeness, we prove the following 

$\  $

{\bf Lemma 2.1.1.} {\it Let $\psi_t$ be a $C^\infty$ topologically transitive Anosov flow, then 
we have the following alternative,\\
$(i)$ $\psi_t$ is topologically mixing,\\ 
$(ii)$ $\psi_t$ admits a $C^\infty$ closed global section with constant return time.}

{\it Proof.} If the case $(ii)$ is 
true, then up to a constant change of time scale, $\psi_t$ is $C^\infty$ flow equivalent to 
the suspension of an Anosov diffeomorphism. Thus it is not topologically mixing. So the alternative is exclusive.

If $\exists \  x\in M$, such that ${W^+_x}$ is not dense in $M$, then 
by Theorem $1.8.$ of {\bf [Pl]}, $E^+\oplus E^-$ is 
integrable with compact $C^1$ leaves. Then the 
{\it canonical $1$-form} of $\psi_t$, $\lambda$, is $C^1$ and in addition 
$$d\lambda =0.$$
So we can find a $C^\infty$ closed $1$-form $\beta$ and a $C^2$ function $f$, such that 
$$\lambda = \beta + df.$$
Denote by $X$ the generator of $\psi_t$, then we have 
$$\lambda(X)-\beta(X) = 1-\beta (X) = X(f).$$
Thus by the cocycle regularity theorem, $f$ is seen to be $C^\infty$ (see {\bf [LMM]}). So $\lambda$ is $C^\infty$. 
We deduce that $E^+\oplus E^-$ is in fact $C^\infty$. Thus if there exists 
$x\in M$ such that $W^+_x$ is not dense in $M$, then the case $(ii)$ of the alternative is true. 

Suppose on the 
contrary that  for $\forall\  x\in M$, $W^+_x$ is dense in $M$. Fix a Riemannian metric on $M$. For $\forall\  
x\in M,\  r >0$, denote by $B(x, r)$ and $B^+(x, r)$ the balls of center $x$ and radius $r$ in $M$ and 
$W^+_x.$ Take arbitrarily two open subsets $U$ and $V$ in $M$ and a small ball $B(y, \epsilon)$ in $V$. Since 
each unstable leaf is dense in $M$ and $M$ is compact, then we can find $R < +\infty$, such that 
$$B^+(x, R) \cap B(y, \epsilon)\not =\emptyset, \  \forall \  x\in M.$$
Take a small disk $B^+(x, \delta)$ in $U$, then by the definition of an Anosov flow, $\exists \  T >0$, such 
that 
$$\psi_t( B^+(x, \delta ))\supseteq B^+(\psi_t(x), R), \  \forall \  t\geq T.$$   
So $\psi_tU\cap V\not =\emptyset, \  \forall \  t\geq T,$ i.e. $\psi_t$ is 
topologically mixing.  $\square$

$\ $

Let us recall the following results established in {\bf [Sa]}.

$\ $

{\bf Theorem 2.1.1.} ({\bf [Sa]}, Theorem $1.3.$) {\it Let $f$ be a topologically 
transitive $C^\infty$ Anosov diffeomorphism ($\phi_t$ be a topologically mixing $C^\infty$ Anosov flow) on a closed 
manifold $M$ which is uniformly quasiconformal on the unstable distribution. Then it is conformal with respect to 
a Riemannian metric on this distribution which is continuous on $M$ and $C^\infty$ along the leaves of the unstable 
foliation.}

$\ $

{\bf Theorem 2.1.2.} ({\bf [Sa]}, Theorem $1.4.$) {\it Let $f$ ($\phi_t$) be a $C^\infty$ Anosov 
diffeomorphism (flow) on a closed manifold $M$ with dim$E^+\geq 2$. Suppose that it is conformal with respect to 
a Riemannian metric on the unstable distribution which is continuous on $M$ and $C^\infty$ along the leaves of 
the unstable foliation. Then the (weak) stable holonomy maps are conformal and the (weak) stable distribution is 
$C^\infty$.}

$\ $     

{\bf Lemma 2.1.2.} {\it Let $\phi_t$ be a $C^\infty$ topologically transitive uniformly quasiconformal Anosov flow, 
such that $E^+\oplus E^-$ is $C^\infty$ and the dimensions of $E^+$ and $E^-$ are at least $2$. Then 
$E^+$ and $E^-$ are both $C^\infty$.}

{\it Proof.} If $\phi_t$ is topologically mixing, then by Theorems $2.1.1.$ and $2.1.2$, 
$E^{+,0}$ and $E^{-,0}$ are both $C^\infty$. Since $E^+\oplus E^-$ is supposed to be $C^\infty$, then 
$E^+$ and $E^-$ are also $C^\infty$. 

If $\phi_t$ is not topologically mixing, then by Lemma $2.1.1,$ $E^+\oplus E^-$ is integrable with 
$C^\infty$ compact leaves. 
Take a leaf of $E^+\oplus E^-$, $\Sigma$, and $T >0$, such that $\phi_T(\Sigma) =\Sigma.$
Then $\phi_T$ is a 
$C^\infty$ topologically transitive uniformly quasiconformal Anosov diffeomorphism. Again by Theorems 
$2.1.1.$ and $2.1.2$, the unstable and stable distributions of $\phi_T$ are $C^\infty$. So 
$E^+$ and $E^-$ of $\phi_t$ are also $C^\infty.$   $\square$

$\ $

{\bf Lemma 2.1.3.} {\it Let $\phi_t$ be as in the lemma above, then $\phi_t$ preserves two $C^\infty$ conformal 
structures along $\mathcal{F}^+$ and $\mathcal{F}^-$, denoted by $\tau^+$ and $\tau^-$.}

{\it Proof.} See {\bf [Ka]} for some details about conformal structures. 
If $\phi_t$ is topologically mixing, then by Theorem $2.1.1$, $\phi_t$ preserves a 
continuous conformal structure $\tau^+$ along $\mathcal{F}^+$, which is $C^\infty$ along the 
leaves of $\mathcal{F}^+$. Let $y\in W^{-, 0}_x$, consider the weak stable holonomy map 
along the weak stable foliation
$$H^{-, 0}_{x, y}: W^+_x \rightarrow W^+_y$$  
$$z\in W^+_x \rightarrow W^+_y\cap W^{-, 0}_z.$$

By Lemma $2.1.2$, $E^+$ and $E^-$ are $C^\infty$. Then $H^{-, 0}_{x, y}$ depends smoothly on 
$x$ and $y$. By Theorem $2.1.2$, $\tau^+$ is invariant under the 
weak stable holonomy maps. So $\tau^+$ is in fact $C^\infty$ on $M$. Similarly, we get a $C^\infty$ 
conformal structure $\tau^-$ along $\mathcal{F}^-$. 

If $\phi_t$ is not topologically mixing, then we get the smooth $\tau^+$ and $\tau^-$ by similar 
arguments about the induced Anosov diffeomorphism on a leaf of $E^+\oplus E^-.$  $\square$

$\ $

Let $\phi_t$ be as in Lemma $2.1.2$. For $\forall\  x\in M$, we extend the conformal structure $\tau^+_x$ 
at $0\in E^+_x$ to all other points of $E^+_x$ via translations. We denote by $\sigma^+_x$ 
this translation-invariant conformal 
structure on $E^+_x$. Then by Lemma $3.1.$ of {\bf [Sa]}, for $\forall \  x\in M$, $h^+_x$ 
sends $\tau^+\mid_{W^+_x}$ to $\sigma^+_x$. So for $\forall \  y\in W^+_x$, $h^+_x\circ ( h^+_y)^{-1}$ is a conformal 
diffeomorphism of $( \sigma^+_y,\  E^+_y)$ to $(\sigma^+_x, \  E^+_x)$. Since the dimension of $E^+$ is at least 
$2$, then $h^+_x\circ (h^+_y)^{-1}$ is naturally an affine map. So if, by 
$h^+_x$ and $h^+_y$, we pull back the canonical flat 
linear connections of $E^+_x$ and $E^+_y$ onto $W^+_x$, we get the same $C^\infty$ connection on $W^+_x$. Thus 
in this way, we get 
a well-defined transversally continuous connection along $\mathcal{F}^+$, denoted by $\nabla^+$ (see 
{\bf [Ka1]} for some details about connections along a foliation).

By the condition $(i)$ of the 
{\it unstable linearization}, $\nabla^+$ is seen to be $\phi_t$-invariant. Similarly, we get a continuous 
$\phi_t$-invariant connection along $\mathcal{F}^-$, $\nabla^-$. If the {\it linearizations}, $\{h^\pm_x\}_{x\in M}$, 
depend smoothly on $x$, then $\nabla^+$ and $\nabla^-$ are certainly $C^\infty$. 
But in general, we can only see that $\{h^\pm_x\}_{x\in M}$ depend continuously on $x$, although 
$E^+$ and $E^-$ are both smooth. \\\\ 
{\bf 2.2. Two dynamical lemmas.}

Let $\phi_t$ be a $C^\infty$-flow on a closed manifold $M$ and $\nu$ be a $\phi_t$-invariant probability 
measure. If 
$\phi_t$ is ergodic with respect $\nu$, then by the Multiplicative Ergodic Theorem of Oseledec, there exists a 
$\nu$-conull $\phi_t$-invariant subset $\Lambda$ of $M$ and a $\phi_t$-invariant measurable 
(Lyapunov) decomposition of $TM\mid_{\Lambda}$,
$$ TM\mid_{\Lambda} = \oplus_{1\leq i\leq k}E_i,$$
such that for $\forall\  u_i\in E_i$, 
$$\lim_{t\to \pm\infty}t^{-1}\log\parallel D\phi_t(u_i)\parallel =\chi_i.$$
Here $E_i$ is called a Lyapunov subbundle and $\chi_i$ its Lyapunov exponent. $E_i$ is also denoted by 
$E_{\chi_i}$. 

By definition, we say that the Lyapunov decomposition of $\phi_t$, with respect to $\nu$, is smooth, if 
there exists a $\nu$-conull $\phi_t$-invariant subset $\Lambda_1$ of $M$ and a $C^\infty$ decomposition of 
$TM$ into smooth subbundles
$$TM =\oplus_{1\leq i\leq k}\bar{E_i},$$
such that the Lyapunov decomposition is defined on $\Lambda_1$ and 
$$\bar{E}_i\mid_{\Lambda_1} = E_i\mid_{\Lambda_1}, \  \forall \  1\leq i\leq k.$$
If in addition the support of $\nu$ is $M$, then this $C^\infty$ decomposition of $TM$ is certainly 
unique and $\phi_t$-invariant. By convention, $\bar{E}_a= E_a: =\{0\}$, if $a$ is not a Lyapunov exponent.

The following lemma is proved in {\bf [BFL1]} (see Lemma $2.5.$ of {\bf [BFL1]}).

$\ $

{\bf Lemma 2.2.1.} {\it Under the notations above, we suppose that the Lyapunov decomposition of 
$\phi_t$ is $C^\infty$ and the support of $\nu$ is $M$. If $K$ is a $C^\infty$ $\phi_t$-invariant 
tensor of type $(1,\  r)$, then we have 
$$K( \bar{E}_{\chi_{i_1}}, \cdots, \bar{E}_{\chi_{i_r}})\subseteq 
\bar{E}_{\chi_{i_1}+\cdots +\chi_{i_r}}.$$ 
Let $\nabla$ be a $\phi_t$-invariant $C^\infty$ connection on $M$, such that 
$$\nabla \bar{E}_i\subseteq \bar{E}_i, \  \forall\  1\leq i\leq k.$$
Then we have $\nabla K =0$, iff 
$\nabla_{\bar{E}_0}K =0$, where 
$E_0$ denotes the Lyapunov subbundle of Lyapunov exponent zero.}

$\ $

If $\phi_t$ is an Anosov flow, then $E_0 =\mathbb{R}X,$ 
where $X$ denotes the generator of $\phi_t$. If 
$\phi_t$ is a contact Anosov flow, i.e. it preserves a $C^\infty$ contact form, then 
by the Anosov property, this contact form must be colinear with the {\it canonical $1$-form} of 
$\phi_t$.

Now we prove the following lemma about general flows.

$\ $

{\bf Lemma 2.2.2.} {\it Let $X$ be a $C^\infty$ vector field on a connected manifold $M$. If 
$f$ is a smooth function on $M$, such that $1 + X(f) >0$, then the flow of $X$ is 
$C^\infty$ flow equivalent to that of $\frac{X}{1+ X(f)}$.}

{\it Proof.} Denote by $\phi_t^X$ the flow of $X$. Thus we can 
construct the following map,
$$ \psi^X : M\longrightarrow M$$
$$ x\to \phi_{f(x)}^X(x).$$
By noting that $X =\frac{\frac{X}{1+X(f)}}{1+(\frac{X}{1+X(f)})(-f)}$, we get a similar 
map $\psi^{\frac{X}{1+X(f)}}$. We observe that
$$D\psi^X(\frac{X}{1+X(f)}) = X.$$
In addition, we can easily verify that 
$$ \psi^X \circ \psi^{\frac{X}{1+X(f)}} = \psi^{\frac{X}{1+X(f)}}\circ  \psi^X =Id_M.$$
So the flow of $X$ is $C^\infty$ flow equivalent to that of $\frac{X}{1+ X(f)}$.  $\square$ \\\\
{\bf 2.3. A geometric lemma}

Let $M$ and $F$ be two $C^\infty$ manifolds. Suppose that the dimension of $M$ is $n$. Denote by 
$\mathcal{F}^{1}(M)$ the frame bundle of $M$. If the general linear group, $GL(n,\  \mathbb R)$, acts 
from left smoothly on $F$, then we get an associated fiber bundle, $\mathcal{F}^1M\rtimes F$ (see {\bf [K-No]} 
for some details about fiber bundles). The $C^\infty$ sections of $\mathcal{F}^1M\rtimes F$ are called the 
geometric structures of type $F$ and order $1$ on $M$. Given a $C^\infty$ linear connection $\nabla$ on $M$, 
we get a horizontal distribution on $\mathcal{F}^1M\rtimes F$, denoted by $\mathcal{H}_F$ (see {\bf [K-No]}). 
Then a geometric 
structure $\sigma$ is called $\nabla$-parallel, if $D\sigma(TM)\subseteq \mathcal{H}_F$. 
Let $\gamma$ be a $C^\infty$ curve in $M$, then $\sigma$ is called parallel along $\gamma$, if 
$D\sigma(\dot{\gamma})\subseteq \mathcal{H}_F$. Thus $\sigma$ is parallel iff it is parallel along each 
smooth curve in $M$.    

Let $\phi$ be a $C^\infty$ diffeomorphism of $M$, then $\phi$ acts on $\mathcal{F}^1(M)$ by its differential, $D\phi$. 
For $\forall\  [\alpha,\  a]\in \mathcal{F}^1M\rtimes F$ (see {\bf [K-No]}), 
$$ \phi_{\ast}[\alpha,\  a] := [(D\phi)(\alpha),\  a].$$
Then $\phi$ acts naturally on the geometric structures and a geometric structure, $\sigma$, is 
called $\phi$-invariant, if $\phi_{\ast}\sigma =\sigma.$ Now we prove the following

$\ $

{\bf Lemma 2.3.1.} {\it Let $\nabla$ be a $C^\infty$ complete linear connection on a connected and 
simply-connected manifold $M$. Let $\tau$ be a $C^\infty$ $\nabla$-parallel geometric structure 
of order $1$. If $\nabla T=0$ and $\nabla R =0$, then the group of $C^\infty$ 
$\nabla$-affine transformations of $M$ which 
preserve $\tau$ is a Lie group and acts transitively on $M$.}

{\it Proof.} Note that the $T$ and $R$ above denote respectively the torsion and the curvature 
tensor of $\nabla$. Denote by $G_\tau$ the group mentionned in the lemma. Then $G_\tau$ is 
a closed subgroup of the group of affine transformations. So $G_\tau$ is naturally a Lie group. 
Fix $u_0\in \mathcal{F}^1M$ and denote by $P(u_0)$ the Holonomy subbundle of $u_0$ (see {\bf [K-No]}). 
Denote by $G$ the group of affine diffeomorphisms of $M$ which preserve 
$P(u_0)$. Then by the assumptions, $G$ is naturally a Lie group and acts transitively on $M$ (see {\bf [K-No]}). 
So we need only prove that $g$ preserves $\tau$, $\forall\  g\in G$. 

 Suppose that $\tau$ is of type $F$. Take $g\in G$. For $\forall\  x\in M$, $\exists\  u \in P(u_0)$ and $a\in F$, such 
that 
$[u,\  a] = \tau(x)$. Since $g_{\ast}(u)\in P(u_0)$, then there exists a piecewise smooth horizontal curve, $u(t)$, in 
$P(u_0)$, such that $u(0) =u$ and $u(1) =g_{\ast}(u).$ Project $u(t)$ to M and 
denote the resulting curve by $\gamma$. Then $[ u(t),\  a]$ gives a 
horizontal lift of $\gamma$ in $\mathcal{F}^1M\rtimes F$ and $[ u(0),\  a] =\tau(x).$ Since $\tau$ is 
$\nabla$-parallel, then 
$\tau\circ \gamma$ is also a horizontal lift of $\gamma$. 
We deduce that $[ u(t),\  a] = (\tau\circ\gamma)(t)$, $\forall \  t\in [0,\  1]$, by the 
uniqueness of the horizontal lift beginning at a fixed point. Then we have 
$$ \tau(g(x)) =\tau (\gamma(1)) =[u(1),\  a]$$
$$=[g_{\ast}(u),\  a] = g_{\ast}[u,\  a] = g_{\ast}(\tau(x)).$$
So $G\subseteq G_\tau$. We deduce that $G_{\tau}$ acts transitively on $M$.  $\square$\\\\
{\bf 3. The proof of Theorem 2.}\\
{\bf 3.1. Preparations.}

Let $\phi_t$ be an Anosov flow on a closed manifold $M$, which satisfies the conditions of Theorem $2$. Then 
by Lemma $2.1.2$, $E^+$ and $E^-$ are $C^\infty$. 
Up to finite covers, we suppose that $M$, $E^+$ and $E^-$ are all orientable. So if we denote by $\mu$ the 
Bowen-Margulis measure 
of $\phi_t$, $\mu$ is given by the integration of a nowhere-vanishing $C^\infty$ volume form on 
$M$ (see {\bf [L-S]} and {\bf [So]}). 

At first, we suppose that $\phi_t$ is topologically mixing. Then G. A. Margulis proved (see {\bf [M]}) that there 
exist two (unique up to scalars) families of measures, $\mu^{\pm, 0}$, supported by the leaves of 
$\mathcal{F}^{\pm,0}$, such that 
$$\mu^{\pm, 0}\circ \phi_t =e^{\pm ht}\mu^{\pm, 0},$$
where $h$ denotes the topological entropy of $\phi_t$. For $\forall\  x\in M$ and $y\in W^-_x$, consider the 
stable holonomy map along the stable foliation 
$$H^-_{x, y}: W^{+, 0}_x\to W^{+, 0}_y$$
$$z\in W^{+, 0}_x\to W^{+, 0}_y\cap W^-_z.$$
By {\bf [M]}, $\mu^{+,0}$ is invariant under the stable holonomy maps. Similarly, $\mu^{-, 0}$ is 
invariant under the unstable holonomy maps.

There exist also two families of measures, $\mu^{\pm}$, supported by the leaves of $\mathcal{F}^\pm$, such 
that 
$$\mu^{\pm}\circ \phi_t =e^{\pm ht}\mu^\pm.$$
In addition, they are absolutely continuous with respect to the weak 
stable and weak unstable holonomy maps (see {\bf[HK]}). Fix a $C^\infty$ 
Riemannian metric on $M$. Denote the induced Riemannian volume forms along $\mathcal{F}^\pm$ 
and $\mathcal{F}^{\pm,0}$ by $\nu^\pm$ and $\nu^{\pm, 0}$. 

$\ $

{\bf Lemma 3.1.1.} {\it Under the notations above, there exist $C^\infty$ positive functions on $M$, $f^\pm$ and 
$f^{\pm, 0}$, such that 
$$\mu^\pm =f^\pm \nu^\pm,\  \mu^{\pm, 0} =f^{\pm, 0}\nu^{\pm, 0}.$$}

{\it Proof.} Let us prove at first that along each leaf of 
$\mathcal{F}^+$, $\mu^+$ is absolutely continuous with respect to $\nu^+,$ 
denoted by $\mu^+\ll\nu^+$. 

Take a small ball 
$B^+(x,\  \delta)$ in $W^+_x$ and $A\subseteq B^+(x,\  \delta)$, such that 
$$\nu^+(A) =0.$$ 
By Lemma $2.1.2$, $E^+$ and $E^-$ are $C^\infty$. Then if $\delta$ is sufficiently small, the following map is 
a well-defined local diffeomorphism for $\forall\  x\in M$ :
$$ B^{-, 0}(x,\  \delta)\times B^+(x,\  \delta)\to M,$$
$$(y,\  z)\to B^+(y,\  2\delta)\cap B^{-, 0}(z,\  2\delta).$$  

Set $\Omega : = B^{-, 0}(x,\  \delta)\times A$, via the local 
diffeomorphisms above. Since $E^+$ is $C^\infty$, then 
for $\forall\  y\in B^{-, 0}(x,\  \delta)$, $\nu^+(y\times A) =0.$ Then by the Theorem of Fubini, $\Omega$ is of 
Lebesgue measure zero. Thus $\mu(A) =0$ by the 
assumption. But $\mu$ can be viewed as a product of $\mu^+$ and $\mu^{-, 0}$ 
(see {\bf[Fo]}), then 
$$\mu^+(y\times A)=0, \  \mu^{-,0}-a.e.$$
Since $\mu^+$ is absolutely continuous 
with respect to the weak stable holonomy maps, then $\mu^+(A) =0$. So $\mu^+\ll\nu^+.$ 
Similarly, we have $\nu^+\ll\mu^+.$ Thus we can find a 
measurable function $f^+$ on $M$, such that $f^+ >0$ and 
$$\mu^+=f^+ \nu^+,$$
i.e. for $\forall\  x\in M$, $\mu^+ = f\mid_{W^+_x}\nu^+.$

For $\forall \  (x,\  t)\in M\times \mathbb{R}$, define 
$$f(x,\  t) := \frac{d\nu^+}{d(\nu^+\circ\phi_t)}(x).$$
Then $f$ is easily seen to be a multiplicative cocycle, positive and 
$C^\infty$ (see {\bf [Fo]} for the definition of a cocycle). We have 
$$f^+\circ \phi_t =\frac{d(\mu^+\circ\phi_t)}{d(\nu^+\circ\phi_t)}$$
$$=\frac{d(\mu^+\circ \phi_t)}{d\mu^+}\cdot \frac{d\mu^+}{d\nu^+}\cdot \frac{d\nu^+}{d(\nu^+\circ\phi_t)} 
= e^{ht}\cdot f^+\cdot f(\cdot,\  t).$$

Then the Livsic cohomological theorem shows that $f^+$ can be taken to be continuous and positive (see 
{\bf [Liv]} and {\bf [Fo]}). Since $f$ is $C^\infty$, then by the cocycle regularity theorem, $f^+$ 
can be taken to be $C^\infty$ (see {\bf [LMM]} and {\bf [Wal]}). Similarly we get the smooth and positive 
functions, $f^-$ and $f^{\pm, 0}$.  $\square$

$\ $

{\bf Remark 3.1.1.} The argument above was originally used by P. Foulon in the case of dimension $3$ (see {\bf [Fo]}). 
The product-decomposition above of the Bowen-Margulis measure is common for general Gibbs measures (see 
{\bf [Ha]}).

$\ $

{\bf Lemma 3.1.2.} {\it Under the notations above, $h^+_x$ and $h^-_x$ depend smoothly on $x$. Then in particular, 
$\nabla^+$ and $\nabla^-$ are $C^\infty$ on $M$.}

{\it Proof.} Suppose that dim$E^+ = n\  (\geq 2).$ By Lemmas $2.1.2.$ and $2.1.3,$ $E^+$ is $C^\infty$ 
and $\phi_t$ preserves 
a $C^\infty$ conformal structure $\tau^+$ along $\mathcal{F}^+.$ By the previous lemma, $\mu^+$ is given 
by a family of $C^\infty$ volume forms along $\mathcal{F}^+$. The volume of a frame of $E^+$ is by definition the 
evaluation of $\mu^+$ on this frame. Then by claiming the $\tau^+$-conformal frames of volume $1$ to be 
orthonormal, we get a well-defined $C^\infty$ Riemannian metric along $\mathcal{F}^+$, denoted by $g^+$. 

Take the leafwise Levi-Civita connections of $g^+$, denoted by $\bar{\nabla}^+$. Since $\tau^+$ is 
$\phi_t$-invariant and 
$$\mu^+\circ\phi_t =e^{ht}\mu^+,$$
then 
$$\phi_t^\ast g^+=e^{\frac{2h}{n}t}g^+.$$
We deduce that $\bar{\nabla}^+$ is $\phi_t$-invariant. For $\forall 
\  x\in M$, define 
$$ \bar{h}^+_x : E^+_x\rightarrow W^+_x,$$
$$u\rightarrow exp^{\bar{\nabla}^+}(u).$$
Because of the $\phi_t$-invariance of $\bar{\nabla}^+$, we get 
$$ \bar{h}^+_{\phi_t(x)}\circ D\phi_t =\phi_t\circ \bar{h}^+_x, \  \forall\  x\in M,\  \forall\  t\in \mathbb{R}.$$
Evidently, $\bar{h}^+_x(0)=x$, $D_x(\bar{h}^+_x) = Id$ and $\bar{h}^+_x$ depends smoothly on $x$.  

Since $g^+$ is evidently complete along each leaf of $\mathcal{F}^+$, 
then $\bar{h}^+_x$ is surjective. Fix a Riemannian metric on $M$. Then by the compactness of $M$, 
there exists $\epsilon >0$, such that for $\forall \  
x\in M$, $\bar{h}^+_x\mid_{\{ u\in E^+_x\  \mid \  \parallel u\parallel <\epsilon \}}$ is a $C^\infty$ diffeomorphism 
onto its image. If $t\gg 1$, then $\phi_{-t}$ contracts $E^+$ exponentially. Since we have in addition 
$$\bar{h}^+_x = \phi_t\circ \bar{h}^+_{\phi_{-t}(x)}\circ D\phi_{-t}, \  \forall\  t >0,$$
then $\bar{h}^+_x$ is in fact injective and nowhere singular. We deduce that  
for $\forall \  x\in M$, $\bar{h}^+_x$ is a $C^\infty$ diffeomorphism. 
By the uniqueness of the {\it unstable linearization} 
(see Proposition $2.1.1$), we get 
$$h^+_x =(\bar{h}^+_x)^{-1},\  \forall\  x\in M.$$
So $h^+_x$ depends also smoothly on $x$. Similarly, we get the $C^\infty$ dependence of $h^-_x$ on $x$. Then 
we deduce that $\nabla^+$ and $\nabla^-$ are $C^\infty$ on $M$ (see Subsection $2.1$).  $\square$

$\ $ 

{\bf Remark 3.1.2.} If $\phi_t$ is not topologically mixing and satisfies the 
conditions in Theorem $2$, then by Lemma $2.1.1$, 
$E^+\oplus E^-$ is integrable with $C^\infty$ compact leaves. Take a leaf $\Sigma$ of the foliation 
of $E^+\oplus E^-$ and $T >0$, such that 
$\phi_T(\Sigma )=\Sigma$. Then 
$\phi_T$ is a $C^\infty$ topologically transitive Anosov diffeomorphism on $\Sigma$. In addition 
by Lemma $2.1.2$, the stable and unstable distributions of $\phi_T$ are both $C^\infty$.   

After some evident modifications, Lemmas $3.1.1.$ and $3.1.2.$ are also 
valid for $(\phi_T,\  \Sigma).$ Just as in the case of flow, we get two $C^\infty$ $\phi_T$-invariant connections 
along $\mathcal{F}^\pm_{\Sigma}$, denoted by $\nabla^\pm_{\Sigma}$ (see 
{\bf [K-Sa]} and Subsection $2.1.$). Then there exists on $\Sigma$ a unique 
$C^\infty$ $\phi_T$-invariant connection, $\nabla$, such that for arbitrary $C^\infty$ sections $Y^\pm$ and $Z^\pm$ of 
$E^\pm$, 
$$\nabla_{ Y^\pm}Y^\mp = P^\mp [Y^\pm,\  Y^\mp],\  \nabla_{Y^\pm}{Z^\pm} = (\nabla^\pm_\Sigma)_{Y^\pm}{Z^\pm},$$
where $P^\pm_{\Sigma}$ denote the projections of $T\Sigma$ onto $E^\pm_{\Sigma}$. 

Then by {\bf [BL]}, $\phi_T$ is $C^\infty$-conjugate to a hyperbolic infranilautomorphism (see also {\bf [K-Sa]}). 
Since $\phi_T$ is in addition uniformly quasiconformal, then up to finite covers, $\phi_T$ is $C^\infty$-conjugate 
to a hyperbolic automorphism of a torus. So Theorem $2.$ is true if $\phi_t$ is not topologically mixing 
and satisfies the conditions of Theorem $2$. \\\\
{\bf 3.2. Homogeneity}

Suppose that $\phi_t$ satisfies the conditions of Theorem $2$. Because of the remark above, we suppose 
in addition that $\phi_t$ is topologically mixing. Since $\phi_t$ is uniformly quasiconformal, 
then with respect to $\mu$, its Lyapunov 
exponents are $\{ a^-,\  0,\  a^+\}$, $a^- < 0< a^+$. In particular, the Lyapunov decomposition 
of $\phi_t$ is $C^\infty$. 
Construct a $C^\infty$ connection $\nabla$ on $M$, such that for arbitrary $C^\infty$ sections 
$Y^\pm$ and $Z^\pm$ of $E^\pm$,
$$\nabla X =0,\  \nabla E^\pm\subseteq E^\pm,$$
$$\nabla_{Y^\pm}Z^\mp = P^\mp [Y^\pm,\  Z^\mp],\  \nabla_{Y^\pm}{Z^\pm} = (\nabla^\pm)_{Y^\pm}Z^\pm,$$
$$\nabla_XY^\pm =[X, \  Y^\pm ]+ a^\pm Y^\pm,$$
where $P^\pm$ denote the projections of $TM$ onto $E^\pm$ with respect to the Anosov splitting. Then by a direct 
verification, $\nabla$ is uniquely determined and $\phi_t$-invariant.

We suppose that dim$E^+ = n$ and dim$E^- =m$. Set $\tau:=(X,\  E^\pm,\  \tau^\pm)$. 
Then by Subsection $2.1$, $\tau$ is a $C^\infty$ 
$\phi_t$-invariant geometric structure of order $1$ on $M$.

$\ $

{\bf Lemma 3.2.1.} {\it Under the notations above, $\tau$ is $\nabla$-parallel.}

{\it Proof.} $\tau$ is in a natural sense the sum of the geometric structures, $X$, $(E^+,\  \tau^+)$, $(E^-,\  \tau^-)$. 
Then $\tau$ is $\nabla$-parallel, iff these structures are parallel respectively. Since $\nabla X =0$, 
then $X$ is $\nabla$-parallel.

Consider the structure $(E^+,\  \tau^+)$, denoted by $\sigma^+$. It is easily seen that 
$\sigma^+$ is $\nabla$-parallel, iff $\nabla E^+\subseteq E^+$ and the $\tau^+$-conformal frames are preserved 
by the parallel transport of $E^+$ along each piecewise smooth curve of $M$ (see 
{\bf [K-No]}).

By the definition of $\nabla$, 
the parallel transport of $E^+$ along the orbits of $\phi_t$ is given by 
$$u^+\rightarrow e^{-a^+ t}\cdot D\phi_t(u^+).$$
Since $\tau^+$ is $\phi_t$-invariant, then the $\tau^+$-conformal frames are preserved along the orbits, i.e. 
$\sigma^+\circ \phi_t$ is horizontal. So $D\sigma^+(X)\subseteq \mathcal{H},$ where $\mathcal{H}$ denotes 
the $\nabla$-horizontal distribution of the corresponding fiber bundle.

Take a smooth curve $\gamma$, tangent to $E^+$. The restriction of $\nabla$ to the 
leaves of $\mathcal{F}^+$ is $\nabla^+$. On the leaf containing $\gamma$, $\nabla^+$ is 
equivalent to the canonical flat connection of a vector space. So the $\tau^+$-conformal frames are 
certainly preserved by the parallel transport along $\gamma,$ i.e. 
$\sigma^+\circ \gamma$ is horizontal. Thus $D\sigma^+(E^+)\subseteq \mathcal{H}.$

Take another smooth curve $\gamma$, tangent to $E^-$. Then by the definition of $\nabla$, 
the parallel transport 
along $\gamma$ is given by the differentials of the weak stable holonomy maps. Since $\tau^+$ is invariant 
with respect to these maps (see Subsection 
$3.2.$ of {\bf [Sa]}), then we deduce that $D\sigma^+(E^-)\subseteq \mathcal{H}$. 

So we get $D\sigma^+(TM) = D\sigma^+ (E^+
\oplus E^-\oplus 
\mathbb{R}X) \subseteq \mathcal{H},$ i.e. $(E^+,\  \tau^+)$ is $\nabla$-parallel. Similarly, 
$(E^-,\  \tau^-)$ is also parallel. we deduce that $\tau$ is $\nabla$-parallel.  $\square$

$\ $

By Lemma $3.1.1$, $\mu^{+, 0}$ is given by a $C^\infty$ family of nowhere-vanishing volume forms along 
$\mathcal{F}^{+, 0}$. Thus $\mu^{+,0}$ can be viewed as a 
$C^\infty$ nowhere-vanishing section of $\wedge^{n+1}(E^{+, 0})^{\ast}$.

By claiming that $\mu^{+, 0}(E^-,\cdots ) :=0,$ 
$\mu^{+, 0}$ can also be viewed as a $C^\infty$ $(n+1)$-form on $M$. In any case, $\mu^{+, 0}$ is a 
$C^\infty$ geometric structure of order $1$ and in the following, we switch from one of these viewpoints to another 
without further precision. 

$\ $

{\bf Lemma 3.2.2.} {\it Under the notations above, $\mu^{+,0}$ is $\nabla$-parallel.}

{\it Proof.} Since the only positive Lyapunov exponent is $a^+$ and dim$E^+ =n$, then by the entropy formula of 
Y. Pesin (see {\bf [Ma]}), $h_{top}(\phi_t) =na^+.$ Then we have 
$$ \mu^{+, 0}\circ\phi_t = e^{n a^+ t}\mu^{+, 0}.$$

In fact, $\mu^{+, 0}$ is $\nabla$-parallel, iff it is parallel along all the curves, tangent to $\mathbb{R}X$ or 
$E^+$ or $E^-$ (see the proof of Lemma $3.2.1.$). 

By the definition of $\nabla$, the parallel transport of $E^+$ along an orbit of $\phi_t$ is given by 
$$u^+\rightarrow e^{-a^+ t}\cdot D\phi_t(u^+).$$ 
Since $\nabla X =0$, then $X\circ \phi_t$ is parallel. Fix $x\in M$ and a basis $\{u_i^+\}_{1\leq i\leq n}$ of 
$E^+_x$. Define 
$$S(t) := (X\circ\phi_t)\wedge e^{-a^+ t}D\phi_t(u^+_1)\wedge\cdots\wedge e^{-a^+ t}D\phi_t(u^+_n).$$
Then $S(t)$ is parallel along this $\phi_t$-orbit, i.e. 
$\nabla_{\partial_{t}}S(t) =0.$ Since 
$$\mu^{+, 0}\circ\phi_t = e^{na^+ t}\mu^{+, 0},$$ 
then 
$$(\phi_t)^\ast\mu^{+, 0} = e^{na^+t}\mu^{+, 0},$$
i.e.
$$ \mu^{+, 0}(\phi_t(x)) = e^{na^+ t}\cdot (\phi_{-t})^\ast(\mu^{+, 0}(x)).$$
So with respect to the natural pairing of $\mu^{+, 0}$ and $S$, we have
$$0 = \partial_{t}<\mu^{+, 0}(\phi_t(x)),\  S(t)>$$
$$=<\nabla_{\partial_t}\mu^{+, 0},\  S> +<\mu^{+, 0},\  \nabla_{\partial_t}S>$$
$$=<\nabla_{\partial_t}\mu^{+, 0}, \  S>.$$
So $\nabla_{\partial_t}\mu^{+, 0} =0,$ i.e. $\mu^{+, 0}$ is parallel along the orbits of $\phi_t$. 

Take a smooth curve $\gamma$, tangent to $E^-$ and beginning at $x$. 
Since $\mu^{+,0}$ is invariant under the stable holonomy maps, then we get 
$$\mu^{+,0}(\gamma(t)) =(H^-_{x,\gamma(t)})_\ast(\mu^{+,0}(x)).$$
By the definition of $\nabla$, for $\forall\  u^+\in E^+_x$, the parallel transport of $u^+$ along $\gamma$ 
is obtained by the differentials of the weak stable holonomy maps,
$$ u^+\rightarrow (DH^{-, 0}_{x,\gamma(t)})(u^+).$$

Take a small curve $l$, tangent to $E^+$ and with $u^+$ as the tangent vector at $0$. Fix 
$t$, then for $\forall\ s\ll1$, 
$H^-_{x,\gamma(t)}(l(s))$ and $H^{-,0}_{x, \gamma(t)}(l(s))$ are contained in $W^{+, 0}_{\gamma(t)}\cap 
W^{-, 0}_{l(s)}$. So for $\epsilon\ll 1$, we can find a smooth function $b : [0,\  \epsilon]\rightarrow\mathbb{R}$, such 
that $b(0) =0$ and  
$$ H^-_{x,\gamma(t)}(l(s)) =\phi_{b(s)}(H^{-, 0}_{x,\gamma(t)}(l(s))),\  \forall \  s\in [0,\  \epsilon].$$
By differentiating the relation above with respect to $s$ at $0$, we get a number $a(t)$, such that 
$$DH^-_{x,\gamma(t)}(u^+) = DH^{-,0}_{x, \gamma(t)}(u^+) + a(t)\cdot X(\gamma(t)).$$

Take a basis of $E^+_x$ as above, $\{u^+_1,\cdots, u^+_n\}$. Define 
$$S(t) :=(DH^-_{x,\gamma(t)})(X_x\wedge u^+_1\wedge\cdots\wedge u^+_n).$$
Then by the relation above, we have 
$$S(t) = X(\gamma(t))\wedge (DH^{-,0}_{x,\gamma(t)})(u^+_1)\wedge\cdots\wedge (DH^{-,0}_{x,\gamma(t)})(u^+_n).$$
We deduce that 
$$\nabla_{\partial_t}S =0.$$
Then as above, we have
$$0= \partial_t<\mu^{+,0}(t),\  S(t)> =<\nabla_{\partial_t}\mu^{+, 0},\  S(t)>.$$
So $\mu^{+, 0}$ is parallel along $\gamma$. 

Now take a curve $\gamma$, tangent to $E^+$ and beginning at $x$. On the line 
bundle $\wedge^{n+1}(E^{+,0})^\ast$, 
$\nabla$ induces naturally a connection $\nabla_1$. Then certainly, $\mu^{+,0}$ is parallel, iff 
$\nabla_1\mu^{+,0} =0.$ 

Along each curve $l$, denote by $P^l_{s_1, s_2}$ the parallel 
transport of $\nabla_1$ from $l(s_1)$ to $l(s_2)$. Denote by $\Omega^{+,0}$ the curvature form of 
$\nabla_1$. Then $\Omega^{+,0}$ is a $\phi_t$-invariant $2$-form on $M$. By the Anosov 
property of $\phi_t$, we get 
$$\Omega^{+,0}(X,\  E^\pm) =0,\  \Omega^{+,0}(E^\pm,\  E^\pm) =0.$$
So the restriction $\nabla_1\mid_{W^{+,0}_x}$ is flat. We 
deduce that if two curves are homotopic with fixed endpoints in $W^{+,0}_x$, then 
their parallel transports are the same. For $\forall\  y\in M$, denote by $O_y$ the $\phi_t$-orbit of $y$. 
Fix $t$, then for $\forall\  s >0,$ we have 
$$c(t):= \cfrac{P^\gamma_{0, t}(\mu^{+,0}(x))}{\mu^{+,0}(\gamma(t))} $$
$$=\cfrac{P^{O_{\gamma(t)}}_{0, -s}\circ P^\gamma_{0,t}(\mu^{+,0}(x))}{P^{O_{\gamma(t)}}_{0, -s}
(\mu^{+,0}(\gamma(t)))}$$
$$=\cfrac{P^{\phi_{-s}\circ\gamma}_{0,\  t}(\mu^{+,0}(\phi_{-s}(x)))}
{\mu^{+,0}(\phi_{-s}(\gamma(t)))},$$
where we have used that $\mu^{+,0}$ is parallel along the orbits of $\phi_t$. If $s\rightarrow +\infty$, then 
the length of 
$\phi_{-s}\circ\gamma$ goes to zero. Thus by the compactness of $M$, 
$c(t)$ goes to $1$, if $s\rightarrow +\infty$. So $c(t)=1$, i.e. $\mu^{+,0}$ is parallel along $\gamma$.

So $\mu^{+, 0}$ is parallel along all the smooth curves tangent to $\mathbb{R}X$ or $E^+$ or $E^-$. 
We deduce that $\mu^{+,0}$ is $\nabla$-parallel.  $\square$

$\ $

View $\mu^{+,0}$ as a $C^\infty$ $(n+1)$-form and define $\omega^+ := i_X\mu^{+,0}$. Since $X$ and $\mu^{+,0}$ 
are $\nabla$-parallel, then $\omega^+$ is also $\nabla$-parallel. Similarly, we get 
a $\nabla$-parallel $m$-form, $\omega^-$. 
 
Set $\sigma := (\tau,\  \omega^+,\  \omega^-)$. Then by Lemmas $3.2.1.$ 
and $3.2.2$, $\sigma$ is a $C^\infty$ $\nabla$-parallel 
geometric structure of 
order $1$ on $M$. Let $\widetilde M$ be the universal covering space of $M$. 
Denote by $\tilde\sigma$ and $\widetilde\nabla$ the 
lifts of $\sigma$ and $\nabla$ to $\widetilde M$. 

$\ $

{\bf Lemma 3.2.3.} {\it The group of $\widetilde\nabla$-affine transformations of $\widetilde M$, which 
preserve $\tilde\sigma$, is a Lie group and acts transitively on $\widetilde M$.}

{\it Proof.} Recall that the Lyapunov decomposition of $\phi_t$, with respect to $\mu$, is $C^\infty$. 
Let us prove at first that $\nabla T=0$ and $\nabla R=0.$

Suppose that $K$ is a $C^\infty$ $\phi_t$-invariant tensor of type $(1,\  k)$. 
Take arbitrarily the Lyapunov exponents 
$\{\chi_1, \cdots ,\chi_k\}$ and $C^\infty$ vector fields $\{Y_1,\cdots,Y_k\}$, such that 
$$Y_i\subseteq \bar E_{\chi_i},\  \forall\  1\leq i\leq k.$$
By the definition of $\nabla$ and Lemma $2.2.1$, we have 
$$(\nabla_XK)(Y_1,\cdots,Y_k) =\nabla_XK(Y_1,\cdots,Y_k) - 
\sum_{1\leq i\leq k}K(Y_1,\cdots, \nabla_XY_i,\cdots, Y_k)$$
$$=[X, K(Y_1,\cdots,Y_k)]+(\sum_{1\leq i\leq k}\chi_i)K(Y_1,\cdots, Y_k)- K([X, Y_1] +\chi_1 Y_1,\cdots)\cdots$$
$$=[X, K(Y_1,\cdots,Y_k)]-\sum_{1\leq i\leq k}K(Y_1,\cdots,[X, Y_i],\cdots, Y_k)$$
$$=(\mathcal{L}_XK)(Y_1,\cdots,Y_k) =0.$$
So $\nabla_XK=0.$ Then by Lemma $2.2.1$, $\nabla K=0.$ In particular, we get 
$$\nabla T =0,\  \nabla R =0.$$ 

By Lemma $2.2.3.$ of {\bf [BFL2]}, the $\nabla$-geodesics tangent to $E^+$ or $E^-$ are defined 
on $\mathbb{R}$. Thus we get 
the completeness of $\nabla$ by the following proposition established in {\bf [Fa]},

$\ $

{\bf Proposition 3.2.1.} ({\bf [Fa]}, Lemma $A$) {\it Let $\nabla$ be a $C^\infty$ linear 
connection on a connected manifold $M$ of dimension $n$. Let $X_1, \cdots, X_k$ be complete fields on $M$ and 
$E_1,\cdots, E_l$ be smooth distributions on $M$, such that \\
$(1)$. $\nabla X_i=0,\  \forall\  1\leq i\leq k,\  \nabla E_j\subseteq E_j,\  \forall\  1\leq j\leq k,$\\
$(2)$. $TM = \mathbb{R}X_1\oplus\cdots \mathbb{R}X_k\oplus E_1\oplus\cdots\oplus E_l,$\\
$(3)$. $\nabla R=0,\  \nabla T=0,$\\
$(4)$. For $\forall\  1\leq j\leq k$, the geodesics of $\nabla$ tangent to $E_j$ are all defined on $\mathbb{R},$\\
then $\nabla$ is complete.}

$\ $
   
We deduce that $\widetilde\nabla$ is also complete. Then we conclude by Lemma $2.3.1.$  $\square$

$\ $

Denote by $G$ the Lie group in the previous lemma. Then $G$ can be viewed as the symmetry group of 
our dynamical system. Fix a point $x\in \widetilde M$ and denote by $H$ the isotropy subgroup of $x$. 
Then $G/H\cong \widetilde M.$ Denote by $\Gamma$ the fundamental group of $M$, thus $\Gamma$ is 
contained in $G$ as a discrete subgroup. 

By claiming the $\tau^+$-conformal frames of $\omega^+$-volume $1$ to 
be orthonormal, we can construct as in Lemma $3.1.2.$ 
a $C^\infty$ fiber metric on $E^+$, denoted again by $g^+$. Similarly, we construct 
a $C^\infty$ fiber metric $g^-$ on $E^-$. Then we get a $C^\infty$ Riemannian metric on $M$
$$g :=\lambda^2\oplus g^+\oplus g^-,$$
where $\lambda$ denotes the {\it canonical $1$-form} of $\phi_t$. 
By the definition of $\sigma$, each element of $G$ preserves $\tilde g$. So in a natural way, $G$ is a closed 
subgroup of the isometry group of $\tilde g$. We deduce that $H$ is a compact Lie subgroup of $G$ 
(see {\bf [Be]} ch.I, $1.78$). \\\\
{\bf 3.3. Symmetric Anosov flows.}

In this subsection, we finish the proof of 
Theorem $2.$ The arguments are based on {\bf [To1]} and {\bf [To2]}. Let 
us recall at first the following definitions (see {\bf [To2]}).

{\bf Definition.} Let $\psi_t$ be a $C^\infty$ flow on $M$, then a Lie transformation group $G$ of $M$ is 
called a symmetric group of $(M,\  \psi_t)$, if $G$ centralizes $\{\psi_t\}$ in Diff$(M)$ and the isotropy 
subgroups are compact in $G$.

The flow $\psi_t$ is called {\it symmetric}, if there exists a normal covering space 
$\bar M$ of $M$, such that the group of deck transformations is contained as a discrete subgroup in an effective and 
transitive symmetric group of the lifted flow, $(\bar M,\  \bar{\psi_t}).$

$\ $

Compared to {\bf [To2]}, we have added to the definition 
the effectiveness of the symmetric group action, which makes no essential difference. 

Let $\psi_t$ be a {\it symmetric} flow. Denote 
by $G$ the symmetric group of $\bar\psi_t$ on $\bar M$. 
Fix a point in $\bar M$ and denote its isotropy subgroup by $K$. 
The Lie algebras of $G$ and $K$ are denoted by $\mathfrak{g}$ and $\mathfrak{g}_K$. Denote the deck group 
by $\Gamma$, thus $M\cong \Gamma
\diagdown G/ K.$ Then by Proposition $1.$ of {\bf [To2]}, $\exists \  \alpha 
\in \mathfrak{g}$, such that $[\alpha,\  \mathfrak{g}_K] \equiv 0$ and 
$$\psi_t(\Gamma g K) =\Gamma(g\cdot exp (t\alpha)) K,\  \forall \  t\in \mathbb{R}.$$
Now by the fundamental estimation in Theorem $1$. of {\bf [To1]}, we get easily    

$\ $

{\bf Proposition 3.3.1.} {\it Under the notations above, if $\psi_t$ is Anosov, then the kernel of $ad\alpha$ is 
$\mathfrak{g}_K
\oplus\mathbb{R}\alpha$ and $ad\alpha$ has no nonzero imaginary eigenvalues. If $\psi_t$ is in addition 
uniformly quasiconformal, then $\Re(Spec(ad\alpha))$ has only three elements.}

$\ $

Given two flows $\phi_t^1$ and $\phi_t^2$. They are called {\it commensurable}, if some finite normal 
cover of $\phi_t^1$ is $C^\infty$ flow equivalent to some finite normal cover of $\phi_t^2$. 
By combining Theorems $5.$ and $6.$ of {\bf [To2]}, we get the following 

$\ $

{\bf Proposition 3.3.2.} {\it Under the notations above, if $\psi_t$ is neither a 
suspension nor a contact flow, then up to commensurability, its lift $\bar{\psi}_t$ can be constructed as follows.

Define $G := N\rtimes(Spin(n, 1)\times K_1\times\cdots\times K_p)$, where $N$ is a vector group 
of positive dimension and $K_1,\cdots, K_p$ are compact, simply connected and almost simple Lie groups. Let 
$\mathfrak{k}'\oplus \mathfrak{p}$ be a Cartan decomposition of $\mathfrak{so}(n,\  1)$ and $\alpha$ be a non-zero 
element of $\mathfrak{p}$. Let $\mathfrak{k}$ be the centralizer of $\alpha$ in $\mathfrak{k}'$ and $K$ be the 
connected Lie subgroup of $G$ with Lie algebra $\mathfrak{k}\oplus\mathfrak{g}_{K_1}\oplus\cdots \oplus 
\mathfrak{g}_{K_p}$. Then we have $Ker(ad\alpha) = \mathfrak{g}_K\oplus \mathbb{R}\alpha$ and  
$$\bar{\psi_t}(g K) =(g\cdot exp (t\alpha)) K,\  \forall\  t\in \mathbb{R},\  \forall\  g\in G.$$}

Now suppose that $\phi_t$ satisfies the conditions of Theorem $2.$ and is topologically mixing. Then by 
Subsection $3.2$, $\phi_t$ turns out to be a {\it symmetric} Anosov flow. 

By Proposition $3.$ of {\bf [To2]} and the previous proposition, it is 
easily seen that up to finite covers, each contact {\it symmetric} Anosov flow 
must be $C^\infty$ flow equivalent to the geodesic flow of a locally symmetric 
space of rank $1$. So if $\phi_t$ is contact, it is 
finitely covered by the geodesic flow of a locally symmetric 
space of rank $1$. Since $\phi_t$ is in addition uniformly quasiconformal, then the locally symmetric space in 
question must have constant negative curvature. Now we finish the proof of Theorem $2.$ by proving 

$\ $

{\bf Lemma 3.3.1.} {\it Suppose that $\phi_t$ satisfies the conditions of Theorem $2$. 
If in addition it is topologically mixing, then $\phi_t$ must be contact.}

{\it Proof.} Suppose on the contrary that $\phi_t$ is not contact. Since 
$\phi_t$ is topologically mixing, then it is not a suspension either. Thus using the notations of Proposition 
$3.3.2$, up to commensurability, a lift of $\phi_t$ is given by 
$$G/K\stackrel{\bar{\phi_t}}{\rightarrow} G/ K $$
$$ gK\to (g\cdot exp t\alpha) K,$$
such that $Ker(ad\alpha) = \mathfrak{g}_K\oplus \mathbb{R}\alpha$.

Since $\alpha\in \mathfrak{p}$ and $\mathfrak{so}(n,1)$ is of rank $1$, then there exists $a>0$ and 
$X_{\pm}\in\mathfrak{so}(n, 1)$, such that 
$$[a\cdot\alpha,\  X_{\pm}] =\pm 2 X_{\pm},\  [X_+,\  X_-] = -a\cdot\alpha.$$

Denote $a\cdot\alpha$ again by $\alpha$, i.e. consider the flow given by $a\cdot\alpha$. 
Denote by $\mathfrak{g}_{\alpha}$ the Lie subalgebra generated by $\{\alpha,\  X_+,\  X_-\}$, then 
$$\mathfrak{g}_{\alpha}\cong \mathfrak{sl}(2,\mathbb{R}).$$
Recall that $G = N\rtimes(Spin(n ,1)\times K_1\times\cdots\times K_p),$ where $N$ is a vector group of 
positive dimension. 
By identifying $N$ with its Lie algebra, we get from this semidirect product a linear representation of 
$\mathfrak{so}(n, 1)$ on $N$. The restriction onto $\mathfrak{g}_{\alpha}$ of 
this representation gives a $\mathfrak{sl}(2,\mathbb{R})$-module (see {\bf [Bo1]} $3$). 

Now by Proposition $2.$ of ({\bf [Bo2]} ch.VIII, $1.2$), there exists a nonzero vector $e$ in $N$ and 
$m\in \mathbb{Z}^+\cup \{ 0\}$, such that 
$$ \alpha(e) = m\cdot e.$$
Since $Ker(ad\alpha) = \mathfrak{g}_K\oplus\mathbb{R}\alpha,$ then $Ker(ad\alpha)\cap N
=\{0\}$ (see Proposition $3.3.2$). So $m\not =0.$ 
By Proposition $3.3.1$, $\Re(spec(ad\alpha))$ has only three elements. 
In addition, $X_+$ is taken such that 
$$[\alpha,\  X_+]=2X_+.$$
We deduce that $m =2.$ Define $e_1=-X_-(e)$, 
then by Propositions $1.$ and $2.$ of ({\bf [Bo2]} ch.VIII, $1.2$), we get 
$$ e_1\not =0,\  \  \  \alpha(e_1) =0,$$
which contradicts to $Ker(ad\alpha)\cap N =\{0\}$.  $\square$

$\ $

{\bf Remark 3.3.1.} It is easily seen that each {\it symmetric} Anosov flow preserves a volume form and 
its Bowen-Margulis measure is in the Lebesgue measure class (see {\bf [Bow]}).\\\\
{\bf 4. Proof of Theorem 1.}\\  
{\bf 4.1. A time change.}

In this subsection, we reduce Theorem $1.$ to the case of Theorem $2.$ The essential point is the following 

$\ $

{\bf Lemma 4.1.1.} {\it Let $\phi_t$ be a $C^\infty$ volume-preserving Anosov flow 
on $M$ with smooth distributions, i.e. 
$E^+$ and $E^-$ are both $C^\infty$. Then there exists a smooth time change of $\phi_t$, which has also 
smooth distributions and whose Bowen-Margulis measure is in the Lebesgue measure class.}

{\it Proof.} Fix a $C^\infty$ Lyapunov metric $g$ on $M$. Then $\exists \  b >0$, such that
$$ \parallel D\phi_{\mp t}(u^\pm)\parallel \leq e^{-bt}\parallel u^\pm\parallel, \  
\forall\  t>0,\  \forall \  u^\pm\in E^\pm. \eqno{(\ast)}$$ 

Up to finite covers, we suppose 
that $E^+$ and $E^-$ are both orientable. Denote by $n$ and $m$ 
the dimensions of $E^+$ and $E^-$ and by $\nu^\pm$ the volume 
forms of $g\mid_{E^\pm}$ on $E^\pm$. For $\forall\  x\in M$ and $\forall\  t\in \mathbb{R}$, define
$$ det(D\phi_t)\mid_{E^\pm_x} : =\frac{(\phi_t^\ast\nu^\pm)_{x}}{\nu^\pm_x}.$$
Then by $(\ast)$, we get for $\forall\  t >0,$ 
$$ det(D\phi_t)\mid_{E^+_x} \geq e^{nbt},\   det(D\phi_t)\mid_{E^-_x}\leq e^{-mbt}. \eqno{(\ast\ast)}$$
For $\forall\  x\in M$, we define
$$\phi^\pm(x) := \frac{\partial}{\partial t}\mid_{t=0}log(det(D\phi_t)\mid_{E^\pm_x}).$$
Since $E^\pm$ are both $C^\infty$, then $\phi^\pm$ are both smooth. In addition by $(\ast\ast)$, we get
$$\phi^+\geq nb > 0, \  \phi^-\leq -mb <0.$$

Denote by $X$ the generator of $\phi_t$ and define a 
$C^\infty$ time change $Y :=\frac{X}{\phi^+}$ of $X$. Denote its flow by $\phi^Y_t$. Certainly, 
$\phi^Y_t$ is also a $C^\infty$ volume-preserving Anosov flow. By {\bf [Par]}, the Bowen-Margulis measure of 
$\phi_t^Y$ is in the Lebesgue measure class. Denote by 
$E^-_Y$ and $E^+_Y$ the stable and unstable distributions of $\phi_t^Y.$ 

Take the dual sections $\omega^\pm$ of $\nu^\pm.$ Thus $\omega^+$ and $\omega^-$ are 
nowhere-vanishing $C^\infty$ sections of $\wedge^nE^+$ 
and $\wedge^mE^-$. For each smooth section $Y^-$ of $E^-$, we define the following differential operator 
acting on the sections of $\wedge^nE^+$, 
$$ \nabla_{Y^-} : sec(\wedge^nE^+)\to sec(\wedge^nE^+),$$
$$Y^+_1\wedge\cdots \wedge Y^+_n \to \sum_{ 1\leq i\leq n} Y^+_1\wedge\cdots\wedge P^+[Y^-, Y^+_i]
 \wedge\cdots\wedge Y^+_n.$$

Now we can define a $C^\infty$ $1$-form $\beta^+$ on $M$, such that for arbitrary 
$C^\infty$ sections $Y^\pm$ of $E^\pm$, 
$$ \mathcal{L}_X \omega^+ = \beta^+(X) \omega^+, \  \beta^+(Y^+) =0,$$
$$ (\nabla_{Y^-})\omega^+ = \beta^+(Y^-) \omega^+.$$

Fix $x\in M$ and take a $C^\infty$ diffeomorphism $\psi : \mathbb{R}^n\to W^+_x$. Then we get the following 
smooth map,
$$ \rho: \mathbb{R}\times \mathbb{R}^n\to W^{+, 0}_x,$$
$$ (t, v)\rightarrow \phi_t(\psi(v)).$$  
If $\epsilon\ll 1$, then $\cup_{-\epsilon <t<\epsilon}W^+_{\phi_tx}$ is 
easily seen to be diffeomorphic to $(-\epsilon, \epsilon)\times \mathbb{R}^n$ under $\rho$. In addition, $\rho$ 
sends $\mathcal{F}^+\mid_{\cup_{-\epsilon <t<\epsilon}W^+_{\phi_tx}}$ to the foliation 
$\{t\times\mathbb{R}^n\}_{-\epsilon <t< \epsilon}$ of $(-\epsilon, \epsilon)\times \mathbb{R}^n$.

So on $\cup_{-\epsilon <t<\epsilon}W^+_{\phi_tx}$, we can find a $C^\infty$ connection along the foliation 
$\mathcal{F}^+\mid_{\cup_{-\epsilon <t<\epsilon}
W^+_{\phi_tx}}$, denoted by $\nabla^+_x$. Then we define a $C^\infty$ 
connection $\nabla_x$ on $\cup_{-\epsilon <t<\epsilon}W^+_{\phi_tx}$, such that  
$$ (\nabla_x) X =0,\  (\nabla_x)_{Y^+} Z^+ = (\nabla^+_x)_{Y^+} Z^+,$$
$$ (\nabla_x)_XY^+ = [X, Y^+],$$
where $Y^+$ and $Z^+$ denote arbitrary $C^\infty$ sections 
of $E^+\mid_{\cup_{-\epsilon <t<\epsilon}W^+_{\phi_tx}}$. 
Denote by $\tau_t$ the $\nabla_x$-parallel transport of $E^+_x$ along the $\phi_t$-orbit 
of $x$. Then by the definition of $\nabla_x$, we get 
$$ \tau_t = D\phi_t, \  (\nabla_x)_X\omega^+ = \beta^+(X) \omega^+.$$
Denote by $\Delta_t$ the determinant of $\tau_t$ with respect to $\nu^+$. Then we have 
$$ \Delta_t = det(D\phi_t)\mid_{E^+_x}.$$
By differentiating the two sides of this equality with respect to $t$ at $0$, we get 
$$-\beta^+(X)_x=\phi^+(x).$$
So  
$$ -\beta^+(X) =\phi^+.$$

For each smooth section $Y^-$ of $E^-$, define $\Omega^+(X, Y^-)$ such that
$$(\mathcal{L}_X\circ \nabla_{Y^-} - \nabla_{Y^-}\circ
\mathcal{L}_X - \nabla_{[X, Y^-]})\omega^+ =\Omega^+(X, Y^-)\omega^+.$$
By a direct calculation, we get   
$$\Omega^+(X, Y^-)\circ\phi_{-t} = \Omega^+(D\phi_t X, D\phi_t Y^-),\  \forall\  t\in \mathbb{R}.$$
In addition, we can view $\Omega^+(X,\cdot)$ 
as a smooth section of $(E^-)^\ast$. So by the Anosov property of $\phi_t$, $\Omega^+(X, E^-) \equiv 0.$ 
By a direct calculation, we get
$$ d\beta^+(X, Y^-) =\Omega^+(X, Y^-),\  \forall\  Y^-\subseteq E^-.$$
Thus
$$ d\beta^+ (X, E^-) \equiv 0. \eqno{(\ast\ast\ast)} $$ 
Define 
$$\alpha^+ :=  -\beta^+, \  f := \frac{1}{\alpha^+(X)}.$$
Then we have $\alpha^+(X) =\phi^+$ and $f X =\frac{X}{\phi^+}$. 

By Lemma $1.2.$ of {\bf [LMM]}, the stable distribution of $f X$, $E^-_Y$, is 
given as follows : 
$$ E^-_Y =\{ u^- +\theta(u^-)\cdot X\  \mid\  \forall \  u^-\in E^-\},$$
where $\theta$ is the unique section of $(E^-)^\ast$ which satisfies the following relation,
$$ \mathcal{L}_X(f^{-1}\theta) = f^{-2}\cdot df\mid_{E^-}.$$
Now using $(\ast\ast\ast)$, we can easily verify that $-\frac{\alpha^+}{\alpha^+(X)}\mid_{E^-}$ 
satisfies the relation above about $\theta$. Thus 
$$E^-_Y =\{ u^- -\frac{\alpha^+(u^-)}{\alpha^+(X)}\cdot X\  \mid\  \forall \  u^-\in E^-\}.$$
So $E^-_Y$ is smooth.

Similar to $\beta^+$, we can define a $C^\infty$ $1$-form $\beta^-$. If we denote by $E^+_{\frac{X}{-\phi^-}}$ the 
unstable distribution of $\frac{X}{-\phi^-}$, then by similar arguments as above, we get
$$E_{\frac{X}{-\phi^-}}^+=\{ u^+ - \frac{\beta^-(u^+)}{\beta^-(X)}\cdot X\  \mid \  \forall\  u^+\in E^+ \}.$$
Thus $E_{\frac{X}{-\phi^-}}^+$ is smooth.

Since $\phi_t$ preserves a volume form, then 
$\phi^+$ and $-\phi^-$ are strictly cohomologous, i.e. 
there exists a smooth function $H$, such that
$$ \phi^+ +\phi^- = X(H).$$
Recall that $\phi^+ >0$ and $-\phi^- >0.$ So we have 
$$ \frac{X}{\phi^+} = \frac{X}{-\phi^- + X(H)} = \frac{\frac{X}{-\phi^-}}{1 + (\frac{X}{-\phi^-})(H)}.$$
Since $Y= \frac{X}{\phi^+}$, then by Lemma $2.2.2$, the flow of $\frac{X}{-\phi^-}$ 
is $C^\infty$ flow equivalent to $\phi^Y_t.$ Since $E_{\frac{X}{-\phi^-}}^+$ is smooth, 
then $E^+_Y$ is also smooth.  $\square$
  
$\ $

Suppose that $\phi_t$ satisfies the conditions of Theorem $1$. 
Then by Lemma $2.1.2$, $E^+$ and $E^-$ 
are both $C^\infty$. Denote by $\phi^Y_t$ a smooth time 
change of $\phi_t$ as in Lemma $4.1.1.$ Then $E^+_Y$ and $E^-_Y$ are also $C^\infty$ and the Bowen-Margulis 
measure of $\phi_t^Y$ is in the Lebesgue measure class. 
In addition by {\bf [Sa]}, $\phi_t^Y$ is also uniformly quasiconformal. So 
$\phi_t^Y$ satisfies the conditions of Theorem $2.$ 

Now we finish the proof of Theorem $1.$ as follows. Lift $\phi_t$ 
and $\phi_t^Y$ to the covering where 
$\phi^Y_t$ becomes an algebraic flow. Recall 
that $\phi_t$ is just a time change of $\phi_t^Y$. If we denote by 
$\lambda$ the {\it canonical $1$-form} of $\phi_t$, then
$$ X= \frac{Y}{\lambda(Y)},\  \  \lambda(Y) >0.$$
Since $d\lambda$ is $\phi_t$-invariant, then by the Anosov property, we have 
$$i_Xd\lambda =0.$$
So $i_Yd\lambda =0$. We deduce that
$$\mathcal{L}_Yd\lambda = d(i_Yd\lambda) +i_Yd(d\lambda) =0.$$
Denote in the following by $\lambda_Y$ the {\it canonical $1$-form} of $\phi_t^Y$.
   
If $\phi_t^Y$ is a suspension, then by ${\bf [Fa1]}$, $\lambda$ is closed and its cohomology 
class is propotional to that of $\lambda_Y$. 
So up to a constant change of time scale, $\phi_t$ is $C^\infty$ flow equivalent to $\phi_t^Y$ 
(see Lemma $2.2.2$). Thus Theorem $1.$ is true in this case.

If $\phi_t^Y$ is a geodesic flow, then 
by {\bf [Ham]}, $\exists\  a\in \mathbb{R},$ such that
$$ d\lambda =a\cdot d\lambda_Y.$$
So there exists a $C^\infty$ closed $1$-form $\gamma$, such that 
$$\lambda =a\cdot\lambda_Y +\gamma.$$
If $a\leq 0$, then $\gamma(Y) >0.$ Thus we can find a closed global section of 
$\phi_t^Y$, which is impossible for a contact flow. We deduce that $a >0.$ Since 
$$ X= \frac{Y}{\lambda(Y)}= \frac{1}{a}\cdot\frac{Y}{1+\frac{\gamma}{a}(Y)},$$
then up to a constant change of time scale, $\phi_t$ is $C^\infty$ flow equivalent to a canonical perturbation of 
the geodesic flow of a hyperbolic manifold. So Theorem $1.$ is true.\\\\
{\bf 4.2 Proof of corollary 1. and Theorem 3.}

Let us prove at first Corollary $1.$ stated in the Introduction.

Suppose that $\phi$ satisfies the conditions of Corollary $1.$ 
Denote by $\phi_t$ the suspension of $\phi$. 
Thus $\phi_t$ satisfies the conditions of Theorem $1.$ 
Since $\Sigma$ is a global section of $\phi_t$, then $\phi_t$ can not be the time change of a geodesic flow. 
So by Theorem $1$, up to a constant change of 
time scale, $\phi_t$ is 
finitely covered by the suspension of a hyperbolic automorphism of a torus. 
Thus $\phi$ is finitely covered by a hyperbolic automorphism of a torus.  $\square$

$\ $

To prove Theorem $3,$ we need only prove the following 

$\ $

{\bf Lemma 4.2.1.} {\it Let $\phi_t$ be a $C^\infty$ volume-preserving uniformly quasiconformal 
Anosov flow, such that $E^+$ and $E^-$ are smooth and dim$E^+ =1$, dim$ E^-\geq 2.$ Then 
up to a constant change of time scale and finite covers, $\phi_t$ is $C^\infty$ flow equivalent to the 
suspension of a hyperbolic automorphism of a torus.}

{\it Proof.} By Lemma $4.1.1$, we can find a time change $\phi_t^Y$ whose Bowen-Margulis measure 
is Lebesgue. Since $E^+_Y$ and $E^-_Y$ are also 
$C^\infty$, then by {\bf [Gh1]}, $E^+_Y\oplus E^-_Y$ is integrable with 
smooth compact leaves. Fix a leaf $\Sigma$ of the foliation of $E^+_Y\oplus E^-_Y$ 
and $T\in \mathbb{R}$, such that 
$\phi_T^Y\Sigma =\Sigma.$ Then the Bowen-Margulis measure of $\phi_T^Y$ is 
also in the Lebesgue measure class.

By similar argument as in Lemma $3.1.1$, we can see that the measures 
$\mu^\pm$ are given by $C^\infty$ volume forms along $\mathcal{F}^\pm_{\Sigma}.$ As in Lemma 
$3.1.2$, we get a $C^\infty$ $\phi_T^Y$-invariant connection $\nabla^-_{\Sigma}$ along 
$\mathcal{F}^-_{\Sigma}$ (see also 
Subsection $2.1$). Denote by $Y^+$ the smooth section of $E^+_Y$, such that $\mu^+(Y^+)\equiv 1$. 
Denote by $h$ the topological entropy of $\phi_T^Y$. Since 
$$\mu^+\circ \phi_T^Y = e^h\mu^+,$$ 
then 
$$(\phi_T^Y)_{\ast}Y^+ = e^h Y^+.$$
Thus we get a $C^\infty$ $\phi_T^Y$-invariant connection $\nabla^+_\Sigma$ along 
$\mathcal{F}^+_\Sigma$, such that 
$$ (\nabla^+_\Sigma)_{Y^+}Y^+ =0.$$
So $\phi_T^Y$ preserves a $C^\infty$ connection as in Remark $3.1.2.$ By {\bf [BL]}, $\phi_T^Y$ is 
finitely covered by a hyperbolic automorphism of a torus. Since 
$\phi_t$ is a $C^\infty$ time change with 
smooth distributions of $\phi_t^Y$, then we conclude by Proposition $1$. of {\bf [Fa1]}.  $\square$

$\ $

Now Theorem $3.$ is just a combination of the previous lemma and Theorem $1.$ and the classification 
of three dimensional volume-preserving Anosov flows with smooth distributions in {\bf [Gh]}. 

Finally based on Theorem $1$, we pose the following questions,

$\ $

{\bf Question 1.} {\it Let $\phi$ be a $C^\infty$ volume-preserving uniformly quasiconformal 
Anosov diffeomorphism. Is it rigid if we suppose that dim$E^+ =1$ and dim$E^-\geq 2$ ?}

$\ $

{\bf Question 2.} {\it Let $\phi_t$ be a $C^\infty$ volume-preserving uniformly quasiconformal 
Anosov flow, such that the dimensions of $E^+$ and $E^-$ are at least 2. Does there 
exist a smooth time change of $\phi_t$ which makes $E^+\oplus E^-$ $C^\infty$ ?}\\\\

$\ $   

{\bf Acknowledgements.} The author would like to thank his thesis advisers, P. Foulon and P. Pansu, for the 
discussions and help. He would like also to thank F. Labourie for clarifying discussions and A. Katok for 
mentionning kindly the papers, {\bf [K-Sa]} and {\bf [Sa]}. The auther would like also to thank the referee for many 
helpful comments and suggestions.

$\ $

{\bf References}

$\ $

{\small
{\bf [Be]}, A. L. Besse, Einstein manifolds, {\it Springer, Berlin-Heidelberg-New York,} (1987).
 
{\bf [BFL 1]} Y. Benoist, P. Foulon and F. Labourie, Flots d'Anosov \`a distributions 
de Liapounov diff\'erentiables. I, {\it Ann. Inst. Henri Poincar\'e} 53 (1990) 395-412.

{\bf [BFL 2]} Y. Benoist, P. Foulon and F. Labourie, Flots d'Anosov \`a distributions 
stable et instable diff\'erentiables, {\it J. Amer. Math. Soc.} 5 (1992) 33-74.

{\bf [BL]} Y. Benoist and F. Labourie, Sur les diff\'eomorphismes
d'Anosov affines \`a feuilletages stable et instable diff\'erentiables, 
{\it Invent. Math.} 111 (1993) 285-308.

{\bf [Bo1]} N. Bourbaki, Groupe et alg\`ebre de Lie, ch. 1, {\it Masson, Paris,} 1960.

{\bf [Bo2]} N. Bourbaki, Groupe et alg\`ebre de Lie, ch. 7 and 8, {\it Hermann, Paris}.

{\bf [Bow]} R. Bowen, The equidistribution of closed geodesics, {\it Amer. J. Math.} 94 (1972) 413-423.

{\bf [Fa]} Y. Fang, Geometric Anosov flows of dimension $5$ with smooth distributions, {\it 
preprint of I.R.M.A, No.2003-009, Strasbourg.}

{\bf [Fa1]} Y. Fang, A remark about hyperbolic infranilautomorphisms, 
{\it C. R. Acad. Sci. Paris, Ser. I 336 No.9} (2003) 769-772.

{\bf [Fo]} P. Foulon, Entropy rigidity of Anosov flows in dimension three, {\it Ergod. Th. and Dynam. Sys.} 21 
(2001) 1101-1112.

{\bf [Gh]} \'E. Ghys, Flots d'Anosov dont les feuilletages stables sont diff\'erentiables, {\it 
Ann. Scient. \'Ec. Norm. Sup. (4)} 20 (1987) 251-270.

{\bf [Gh1]} \'E. Ghys, Codimension one Anosov flows and suspensions, {\it Lecture Notes in 
Mathematics} 1331 (1988) 59-72.

{\bf [Ha]} N. T. A. Haydn, Canonical product structure of equilibrium states, {\it Random and 
computational dynamics} 2(1) (1994) 79-96.

{\bf [Ham]} U. Hamenstadt, invariant two-forms for geodesic flows, {\it Math. Ann.} 301 (1995) 
{\it No.4} 677-698.

{\bf [HK]} B. Hasselblatt and A. Katok, Introduction to the modern theory of dynamical 
systems, {\it Encyclopedia of Mathematics and its Applications, vol 54.} 1995. 

{\bf [Ka]} M. Kanai, Differential-geometric studies on dynamics of geodesic and frame flows, {\it 
Japan. J. Math.} 19 (1993) 1-30.

{\bf [Ka1]} M. Kanai, Geodesic flows of negative curved manifolds with smooth stable and unstable foliations, 
{\it Ergod. Th. and Dynam. Sys.} 8 (1988) 215-240.

{\bf [KL]} A. Katok and J. Lewis, Local rigidity for certain groups of toral automorphisms, {\it Israel 
J. Math.} 75 (1991) 203-241.

{\bf [K-No]} S. Kobayashi and K. Nomizu, Foundations of Differential 
Geometry, vol. 1 and 2, {\it Interscience, New York and London,} 1963.

{\bf [K-Sa]} B. Kalinin and V. Sadovskaya, On local and global rigidity of quasiconformal Anosov 
diffeomorphisms, to appear in {\it Journal of the Institute of Mathematics of Jussieu.}

{\bf [Liv]} A. N. Livsic, Cohomology of dynamical systems, {\it Math. USSR Izvestija} 6(6) (1972) 1278-1301.

{\bf [LMM]} R. de La llave, J. Marco and R. Moriyon, Canonical perturbation theory of Anosov systems, and 
regularity results for Livsic cohomology equation, {\it Ann. Math.} 123 (1986) 537-612.

{\bf [L-S]} A. N. Livsic and Ja. G. Sinai, On invariant measures compatible with the smooth structure 
for transitive U-systems, {\it Soviet Math. Dokl.} 13(6) (1972) 1656-1659. 
 
{\bf [M]} G. A. Margulis, Certain measures associated with U-flows on compact manifolds. {\it Func. Anal. Applic.} 
4 (1970) 54-64.

{\bf [Ma]} R. Man\'e, Ergodic theory and differentiable dynamics, {\it Springer Verlag, Berlin, New York,} 1987.

{\bf [Par]} W. Parry, Synchronisation of canonical measures for hyperbolic attractors, {\it Commun. Math. Phys.} 
106 (1986) 267-275.

{\bf [Pl]} J. F. Plante, Anosov flows, {\it Amer. J. Math.} 94 (1972) 729-754.

{\bf [Sa]} V. Sadovskaya, On uniformly quasiconformal Anosov systems, to appear in {\it Math. Res. Lett}.

{\bf [Sh]} M. Shub, Stabilit\'e globale et syst\`emes dynamiques, {\it Ast\'erisque} 56 (1978).

{\bf [So]} V. V. Solodov, Topological topics in dynamical systems theory, {\it Russian Math. Surveys} 46(4) 
(1991) 107-134.

{\bf [Su]} D. Sullivan, On the ergodic theory at infinity of an arbitrary discrete group of hyperbolic motions, in 
``Riemann surfaces and related topics'', {\it Annals of Math. Studies} 97 (1981) 465-497.

{\bf [To1]} P. Tomter, Anosov flows on Infra-homogeneous Spaces, {\it Proc. Symp. in Pure Math, 
Vol. XIV, Global Analysis,} (1970) 299-327.

{\bf [To2]} P. Tomter, On the classification of Anosov flows, {\it Topology} 14 (1975) 179-189.

{\bf [Wal]} C. Walkden, Livsic theorems for hyperbolic flows, {\it Trans. Amer. Math. Soc.} 352 (2000) 1299-1313.

{\bf [Yu]} C. Yue, Quasiconformality in the geodesic flow of negatively curved manifolds. 
{\it GAFA } 6(4) (1996) 740-750.}
\end{document}